\documentclass{article}
\usepackage{amsmath,amstext,amsthm,amsfonts}
\usepackage[dvips]{graphicx}

\theoremstyle{plain}
\newtheorem{theorem}{Theorem}[section]

\newtheorem{lemma}[theorem]{Lemma}
\newtheorem{corollary}[theorem]{Corollary}
\newtheorem{maintheorem}{Theorem}

\theoremstyle{definition}
\newtheorem{remark}[theorem]{Remark}

\newtheorem{definition}[theorem]{Definition}

\newcommand{\field}[1]{\mathbb{#1}}
\newcommand{\real}{\field{R}}

\renewcommand{\natural}{\field{N}}

\newcommand{\SB}{{\cal B}}

\newcommand{\SI}{{\cal I}}

\newcommand{\SK}{{\cal K}}

\newcommand{\SP}{{\cal P}}
\newcommand{\SQ}{{\cal Q}}
\newcommand{\SR}{{\cal R}}

\begin{document}

\title{Equilibrium States for \\ Non-uniformly Expanding Maps}
\author{Krerley Oliveira \footnote{The author is supported by CNPq,
Brazil.}}
\date{\today}

\maketitle

\begin{abstract}
We construct equilibrium states, including measures of maximal
entropy, for a large (open) class of non-uniformly expanding maps
on compact manifolds. Moreover, we study uniqueness of these
equilibrium states, as well as some of their ergodic properties.
\end{abstract}

\section{Introduction}

The theory of equilibrium states originates from statistical
mechanics and was thoroughly developed, in the classical setting
of uniformly hyperbolic dynamical systems, in the seventies and
eighties, especially by Sinai, Ruelle, Bowen, Parry and Walters.

In general, given a continuous transformation $f: M\rightarrow M$
on a compact metric space, and given a continuous function $\phi$,
we call \textbf{equilibrium state} for $(f,\phi)$ a Borel probability
measure $\mu_\phi$ such that
$$
h_{\mu_\phi}(f) + \int \phi d\mu_\phi
=\sup\limits_{\mu \in \SI}\{h_\mu(f) + \int \phi d\mu\},
$$
where the supremum is taken over the set $\SI$ of $f$-invariant
probabilities. That is, an equilibrium state is a maximum of the
function $F_\phi : \SI \rightarrow \field{R}$ defined by
$$
F_\phi(\mu)= h_\mu(f) + \int \phi d\mu.
$$

It is now classical that for uniformly hyperbolic diffeomorphisms,
as well as for uniformly expanding maps, equilibrium states always
exist and they are unique if the potential is H\"older continuous,
assuming that the transformation $f$ is transitive. See
\cite{Pa64,Si72,Bo75,Ru89b}. Moreover, the equilibrium states
coincide with the Gibbs measures, that is, the invariant
probability measures satisfying
\begin{equation} \label{eq1}
\mu(\SB_\epsilon(n,x)) \cong \exp {( \sum\limits_{i=0}^{n-1}
\phi(f^i(x)) - n P)}
\end{equation}
for some $P\in \field{R}$, called the \textit{pressure} of $\phi$,
where $\SB_\epsilon(n,x)$ is the dynamical ball of length $n$ and
size $\epsilon$ around $x$,
$$
\SB_\epsilon(n,x)=\{y\in M; d(f^i(y),f^i(x))\leq \epsilon,
\text{ for every } 0\leq i\leq n-1\},
$$
and $\cong$ means equality up to a uniform factor, independent of
$x$ and $n$. In this setting, the pressure $P$ is given by
$$
P = \sup\limits_{\mu \in \SI} \big\{h_\mu(f) + \int \phi d\mu\big\}.
$$

Several authors have been studying equilibrium states for
non-hyperbolic systems: Bruin, Keller~\cite{BrK98} and Denker,
Urbanski~ \cite{DU92,Ur98}, for interval maps and rational
functions on the sphere, and Buzzi, Maume,
Sarig~\cite{Bu99,BM02,BS,Sa03} and Yuri~\cite{Yu99,Yu00,Yu03}, for
countable Markov shifts and for piecewise expanding maps in one
and higher dimensions, to mention just a few of the most recent
works. Several of these papers, and
particularly \cite{DU92,Sa03,Yu99,Yu00,Yu03},
consider systems with neutral periodic points, a setting of
non-hyperbolic dynamics which has attracted a great deal of
attention over the last years.
Despite all these important contributions, it is fair to say
that the theory of equilibrium states is very much incomplete
outside the uniformly hyperbolic case.

The present work may be seen as a step towards obtaining such a
theory in a general setting of \emph{non-uniformly hyperbolic}
systems (non-zero Lyapunov exponents). Indeed, we prove existence
of equilibrium states for fairly general potentials and for a
robust (open) class of non-uniformly expanding maps. This class
will be defined precisely in the next section. Here we just
mention one of its main features:

\begin{equation} \label{eq2}
\limsup\limits_{n\rightarrow \infty} \frac{1}{n}
\sum\limits_{i=0}^{n-1} \log \|Df(f^i(x))^{-1}\| \leq -2c <0
\end{equation}
for ``most" points, including a full measure set relative to the
equilibrium states $\mu_\phi$ that we construct. We prove that
equilibrium states do exist for every potential $\phi$ with small
variation, that is, such that

\begin{equation}\label{eq3}
\sup \phi-\inf \phi < K
\end{equation}
 for some convenient constant $K$ (see
Definition~\ref{d.low} and comments following it). As a
consequence of \eqref{eq2}, these measures are non-uniformly
expanding, that is, $\mu_\phi$-almost every point has only
positive Lyapunov exponents. Moreover, is possible to prove that
$\mu_\phi$ has a kind of weak Gibbs property (see \cite{Yu00})
as in \eqref{eq1}, with $\cong$ meaning equality up to a factor
with subexponential growth on the orbit of each $x$.

The basic strategy for the construction is to find a
 subset $\SK$ of invariant probability measures which are
expanding and such that $\mu$-almost every point has infinitely
many hyperbolic times, in the sense of \cite{Al00,ABV00}. We prove
that $\nu\mapsto F_\phi(\nu)$ is upper-semicontinuous on $\SK$ and
there exist  maximum on $\SK$ of $F_{\phi}$. Using \eqref{eq3} we
check that the maximum obtained on $\SK$ is really a maximum of
$F_\phi$ over all invariant probabilities.

These arguments apply, in particular, when the potential $\phi$ is
constant, in which case $\mu_\phi$ maximizes the entropy:
$$
h_{\mu_\phi}(f)=h_{top}(f).
$$
However, in this case we can go much further. Using a different
approach, via semi-conjugation to a one-sided subshift of finite
type, we are able to prove that the maximal entropy measure is unique
and a Markov measure. If $f$ is topologically mixing, this measure
is Bernoulli.

Closing this introduction, we mention some questions that are
naturally raised by our results. The first one is to prove
uniqueness or ,at least, finiteness of the equilibrium states for
H\"older potentials in the general situation of \eqref{eq3}, under
topological transitivity. In this direction, in a forthcoming work
the author proves existence and uniqueness of equilibrium states
for an open class of local diffeomorphisms and for potentials with
low variation satisfying a summability condition. Another question
concerns the need of hypothesis \eqref{eq3} itself. Ongoing work
indicates that equilibrium states do exist also for potentials
with large variation, but they may have both positive and negative
Lyapunov exponents: from the viewpoint of these measures the
dynamics looks hyperbolic, rather than expanding.

\section{Setting and  Statements}

We always consider a $C^{1+\alpha}$ local diffeomorphism $f:M^l
\rightarrow M^l$, defined on a compact Riemannian manifold with
dimension $l$. Let $m$ be normalized Lebesgue measure on $M$. We
suppose that $f$ satisfies, for positive constants $\delta_0$,
$\beta$, $\delta_1$, $\sigma_1$, and $p$, $q\in\natural$,

\begin{enumerate}

\item[(H1)] There exists a covering $B_1, \dots,B_p, \dots, B_{p+q}$ of
      $M$ such that every $f|B_i$ is injective and

\begin{itemize}

\item  $f$ is uniformly expanding at every $x\in B_{1} \cup \dots \cup
B_{p}$: $$\|Df(x)^{-1}\| \leq (1+\delta_1)^{-1}.$$

\item  $f$ is never too contracting: $\|Df(x)^{-1}\|\leq (1+\delta_0)$ for
every $x\in
M$.

\end{itemize}

\item[(H2)] $f$ is everywhere volume-expanding:  $|\det Df(x)|\geq \sigma_1$
      with $\sigma_1>q$.

Define $$V=\{x\in M; \|Df(x)^{-1}\|>(1+\delta_1)^{-1}\}.$$

\item[(H3)] There exists a set $W \subset B_{p+1}\cup\dots\cup B_{p+q}$
containing $V$ such that
$$
M_1>m_2 \quad\text{and}\quad m_2-m_1<\beta
$$
where $m_1$ and $m_2$ are the infimum and the supremum
of $\log \|\det Df\|$ on $V$, respectively,
and $M_1$ and $M_2$ are the
infimum and the supremum of $\log \|\det Df\|$ on $W^c$, respectively.
\end{enumerate}

\begin{definition}
 The supremum of $F_\phi$ over the set all invariant probability is called the
\textit{pressure of} $\phi$ and will be denoted by $P(\phi)$
\end{definition}

\begin{definition}\label{d.low} Given $\phi :M \rightarrow \field{R}$
continuous,
we say that $\phi$ has $\rho$-\textit{low variation} if
$$
\max\limits_{x\in M} \phi(x)< P(\phi)-\rho h_{top}(f).
$$
\end{definition}
\begin{remark} 
Note that this is an open condition on the potential, with respect
to the $C^0$ topology. Note also that this is somewhat more general
than condition \eqref{eq3} for $K$ small enough:
assuming $K$ is less than $(1-\rho)h_{top}(f)$, if $\phi$ satisfies
\eqref{eq3} then $\tilde{\phi}=\phi - \inf \phi$ has $\rho$-low
variation potential; the conclusions of Theorem~\ref{theoremA} are
not affected if one replaces $\phi$ by $\tilde{\phi}$, because
potentials that differ by a constant have the same equilibrium states.
\end{remark}

\smallskip

Our first main result is

\begin{maintheorem}\label{theoremA}
Assume hypotheses (H1), (H2), (H3) hold, with $\delta_0$ and
$\beta$ sufficiently small. Then, there exists $\rho$ such that if
$\phi$ is a continuous potential with $\rho$-low variation then
$\phi$ has some equilibrium state. Moreover, these equilibrium
states are hyperbolic measures, with all Lyapunov exponents bigger
than some $c(\delta_1, \sigma_1, p, q)>0$.
\end{maintheorem}

For maximal entropy mesures we are able to say a lot more, under
the following additional hypothesis:

\begin{enumerate}

\item[(H4)] There exists a Markov partition $\SR=\{R_1,\dots,R_d\}$
for $f$ such that
\begin{itemize}
\item $\SR$ is transitive: for any $i,j$ there exists a $k$ such that
$f^k(R_i)\cap R_j \neq \emptyset$;

For simplicity in the proofs we also assume that $W \subset R_1$.
\
\end{itemize}
\end{enumerate}

We say  that a system $(f,\mu)$ is Bernoulli (respectively Markov),
if it is ergodically equivalent to a subshift of finite type endowed
with a Bernoulli (respectively Markov) measure. See e.g. \cite{Man87}
for definitions.

\begin{maintheorem}\label{theoremB}
Assume hypotheses (H1), (H2), (H3), (H4) hold with $\delta_0$ and
$\beta$ sufficiently small. Then there exists a unique invariant
measure $\mu_{max}$ with $h_{\mu_{max}}(f)=h_{top}(f)$. This
measure also satisfies,

\begin{enumerate}

\item all Lyapunov exponents of $\mu_{max}$ are larger than
some $c(\delta_1,\sigma_1,p,q)>0$;

\item $(f,\mu_{max})$ is Markov and, if $f$ is topologically
mixing, it is Bernoulli.

\end{enumerate}
\end{maintheorem}

 Important related results have been obtained by
Yuri~\cite{Yu99,Yu00,Yu03}, where she studies equilibrium states
for very general Markov systems. On the other hand, our hypotheses are
different and, to the best of our knowledge, our approach for
proving Theorem~\ref{theoremA} is new. For one thing, we do not
assume existence of a \emph{generating} Markov partition: instead,
we construct a special partition and \emph{prove} that it is
generating for a carefully chosen class of measures (the set $\SK$
mentioned in the Introduction). In general terms, we exploit the
notion of hyperbolic times to deduce from the dynamical behaviour
certain facts valid at almost every point, uniform versions of
which are taken as hypotheses in Yuri's approach. For instance, we
need no analogue of hypothesis (C5) in \cite{Yu99}: in fact, for
our examples in Section~\ref{s.examples} the diameters of
cylinders \emph{do not} tend to zero. Besides, there a low variation potential
may not satisfies the condition  (C4) in \cite{Yu99}. A combination of both
viewpoints should lead to further progress in this area.

\smallskip

\textbf{Acknowledgements:} I am very thankful to my advisor Marcelo
Viana for his exceptional advice and friendship. Warm thanks go
also to A. Tahzibi, J. Bochi, C. Matheus, and A. Arbieto for
suggestions and many fruitful discussions. I am indebted to IMPA
and its staff for a fine working environment, and to CNPq for
financial support.

\section{Examples}\label{s.examples}

In this section, we sketch the construction of a non-hyperbolic
map $f_0$ that satisfies the conditions in  theorems A and B
above. It will be clear from this construction that these
conditions hold, in fact, for every map $f$ $C^1$-close to $f_0$.

We start by considering any Riemann manifold that supports an expanding map
$g:M \rightarrow M$. For simplicity, choose  $M=\mathbb{T}^n$
  the $n$-dimensional torus, and $g$ an endomorphism induced from a linear
map with eigenvalues $\lambda_n>\dots>\lambda_1>1.$ Denote by $E_i(x)$
the eigenspace
associated to the eigenvalue $\lambda_i$ in $T_x M$.

Since $g$ is expanding, it admits a transitive Markov partition $
R_1, \dots, R_d $ with arbitrarily small diameter.
We may suppose that $g|R_i$ is injective for every $i=1,\dots,d$.
Replacing $g$ by a iterate if necessary, we may suppose that there exists a
fixed point $p_0$ of $g$ and, renumbering if necessary,
this point is contained in the interior of the rectangle $R_d$ of
the Markov partition.

 Considering a small neigbourhood $W \subset R_d$ of $p_0$ we
deform $g$ inside $W$ along the  direction $E_1$. This deformation consists
essentially in rescaling the expansion
along the invariant manifold associated to $E_1$ by a real function $\alpha$.
Let us be more precise:

Considering $W$ small, we may identify $W$ with a neighbourhood of
$0$ in $\real^n$ and $p_0$ with $0$. Without loss of generality, suppose that
$W=(-2\epsilon,2\epsilon)
\times B_{3r}(0)$,
where $B_{3r}(0)$ is the ball or radius $3r$ and
center $0$ in $\real^{n-1}$. Consider a function
$\alpha:(-2\epsilon,2\epsilon)\rightarrow  \real$ such
$\alpha(x)=\lambda_1x$ for every $|x|\geq \epsilon$ and for small constants
$\gamma_1,\gamma_2$:

\begin{enumerate}

\item  $(1+\gamma_1)^{-1}<\alpha'(x)< \lambda_1 + \gamma_2$

\item $\alpha'(x)<1$ for every $x \in
(-\frac{\epsilon}{2},\frac{\epsilon}{2})$;

\item $\alpha$ is $C^0$-close to $\lambda_1$: $\sup\limits_{x\in
(-\epsilon,\epsilon)}|\alpha(x)-\lambda_1x| < \gamma_2$,

\end{enumerate}

Also, we consider a bump function $\theta: B_{3r}(0)\rightarrow \real$ such
$\theta(x)=0$ for every $2r\leq |x| \leq 3r$
 and $\theta(x)=1$ for every $0\leq |x|\leq r$. Suppose that
$\|\theta'(x)\|\leq
C$ for
every $x\in B_{3r}(0)$.  Considering coordinates $(x_1,\dots,x_n)$ such that
$\partial_{x_i} \in E_i$, define
$f_0$ by:

$$ f_0(x_1,\dots,x_n)=
(\lambda_1x_1+\theta(x_2,\dots,x_n)(\alpha(x_1)-\lambda_1x_1),
\lambda_2x_2,\dots,\lambda_n x_n)$$

Observe that by the definition of $\theta$ and $\alpha$ we can extend $f_0$
smoothly to $\field{T}^n$ as $f_0=g$
outside $W$. Now, is not difficult to prove that $f_0$ satisfies the conditions
(H1), (H2), (H3) and (H4) above.

First, we have that $\|Df_0(x)^{-1}\|^{-1}\geq \min\limits_{i=1,\dots,n}
\|\partial_{x_i} f_0\|.$ Observe that:

$$\partial_{x_1} f_0(x_1,\dots,x_n)=
(\alpha'(x_1)\theta(x_2,\dots,x_n)+(1-\theta(x_2,\dots,x_n))\lambda_1,0,\dots,0)$$

$$\partial_{x_i} f_0(x_1,\dots,x_n) =
((\alpha(x_1)-\lambda_1)\partial_{x_i}\theta(x_2,\dots,x_n),0,\dots,\lambda_i,0,\dots,0),
\text{ for } i\geq 2.$$

Then, since $\|\partial_{x_i} \theta(x)\| \leq C$ for every $x\in B_{3r}(0)$,
and $\alpha(x_1)-\lambda_1x_1
\leq \gamma_2$ we have that $ \|\partial_{x_i} f_0\|>(\lambda_i-\gamma_2C)$
for every
$i=2,\dots,n$. Moreover, by condition 1, $\|\partial_{x_1} f_0\| \leq \max
\{\alpha'(x_1), \lambda_1\}\leq
\lambda_1 + \gamma_2,$ if we choose $\gamma_2$ small in such way that
$\lambda_2-\gamma_2C >
\lambda_1 +\gamma_2 $ then:

$$\|\partial_{x_i} f_0\|>\|\partial_{x_1} f_0\|, \text{ for every } i\geq
2.$$

Notice also that $\|\partial_{x_1} f_0\| \geq \min\{\alpha'(x_1),\lambda_1\}
\geq (1+\gamma_1)^{-1}.$ This prove
that:

$$\|Df_0(x)^{-1}\|^{-1}\geq \min\limits_{i=1,\dots,n} \|\partial_{x_i} f_0\| \
(1+\gamma_1)^{-1}.$$ Since $f$ coincides with $g$ outside $W$, we have
$\|Df_0(x)^{-1}\|\leq \lambda_1^{-1}$ for every
$x\in W^c$. Together with the above inequality, this proves  condition (H1),
with $\delta_0=\gamma_1$.

Choosing $\gamma_1$ small and $p=d-1$, $q=1,$ $B_i=R_i$ for every
$i=1,\dots,d$,
 condition (H2) is imediate.
Indeed, observe that the Jacobian of $f_0$  is given by the formula:

 $$\det Df_0(x)=
(\alpha'(x_1)\theta(x_2,\dots,x_n)+(1-\theta(x_2,\dots,x_n))\lambda_1)\prod_{i=2}^n
 \lambda_i.$$  Then,  if we choose $\gamma_1 < \prod_{i=2}^n
 \lambda_i-1$:

 $$\det Df_0(x)>(1+\gamma_1)^{-1}\prod_{i=2}^n
 \lambda_i >1.$$ Therefore, we may take $\sigma_1=
(1+\gamma_1)^{-1}\prod_{i=2}^n \lambda_i>1.$

To verify property (H3) for $f_0$, observe that if we denote by

$$V=\{x\in M; \|Df_0(x)^{-1}\|>(1+\delta_1)^{-1}\},$$

with $\delta_1<\lambda_1-1,$ then $V\subset W.$ Indeed, since $\alpha(x_1)$
is constant equal to
$\lambda_1x_1$ outside $W$  we have that $\|Df_0(x)^{-1}\|\leq \lambda_1^{-1}
< (1+\delta_1)^{-1}$, for every
$x\in W^c$.  Given $\gamma_3$ close to 0, we may choose $\delta_1$ close to 0
and $\alpha$ satisfying the conditions above in
such way that,

$$\sup\limits_{x,y\in V}\alpha'(x_1)-\alpha'(y_1)< \gamma_3.$$ If $m_1$ and
$m_2$ are the
infimum and the supremum of $|\det Df_0|$ on $V$, respectively,

$$m_2-m_1\leq C(\sup\limits_{x,y\in V}
\alpha'(x_1)-\alpha'(y_1)) < \gamma_3 C,$$
where $C=\prod\limits_{i=2}^n\lambda_i$.
Then, we may take $\beta=\gamma_3C$ in (H3). If $M_1$ is the
infimum  of $|\det Df_0|$ on $W^c$, $M_1>m_2$, since $\lambda_1>
(1+\delta_1)\geq \sup\limits_{x\in V}
\alpha'(x)$.  Condition (H4) is clear from the construction, since $f_0=g$
outside $W \subset R_d$, so the Markov
property of $\{R_1,\dots,R_d\}$ is not affected by the pertubation.

    The arguments above show that the hypotheses $(H1), (H2),(H3)$ and $(H4)$
are satisfied by $f_0$. Moreover, if we one takes
    $\alpha(0)=0$, then $p_0$ is fixed point for $f_0$, which is not a
reppeler,
since $\alpha'(0)<1$.
    Therefore, $f_0$ is not a uniformly expanding  map.

    It is not difficult to see that this construction may be carried out in
such
way that $f_0$ does not satisfy the  expansiveness
 property: there is a fixed hyperbolic saddle point $p_0$ such that the stable
manifold of $p_0$ is
 contained in the unstable manifold of two other fixed points.

 For a discussion of related
 examples  see \cite{Ca93} and \cite{BoV00}.

\section{Expanding measures and hyperbolic times}

The proof of Theorem~\ref{theoremA} occupies this section and the
next one. Let us begin by detailing a bit more our strategy to
prove the existence of equilibrium states:

\begin{itemize}

\item To exhibit a  subset $\SK$ of invariant measures such that
all their Lyapunov exponents are positive, and almost every point has
infinitely many hyperbolic times.

\item To show that there exists a common generating partition for
all the measures in $\SK$. This allows us to prove that the
function $\mu \rightarrow h_\mu(f)+\int \phi d\mu$ is
upper-semicontinuous on $\SK$. Using this, we get
that the maximum of $h_\mu(f)+\int \phi d\mu$ for measures $\mu$
in $\SK$ is attained.

\item To prove that if the potential has low variation, the maximum
obtained over $\SK$ is, in fact, a global maximum for the function
$h_\mu(f)+\int \phi d\mu$ over all invariant measures.

\end{itemize}

Let us begin by stating our precise conditions on the constants
$\delta_0$ and $\beta$ in the theorem.
According to \cite[Appendix]{ABV00}, if $f$ satisfies (H1) then
there exists $\gamma_0<1$ depending only on $(\sigma_1,p,q)$ such
that Lebesgue almost every point spends at most a fraction
$\gamma_0$ of time inside $B_{p+1}\cup \dots \cup B_{p+q}$.

Reducing $\delta_0$ if necessary, we may find constants
$\alpha>0$, as close to $1$ as we want, and $c>0$ such that
\begin{equation}\label{eq.alpha}
(1+\delta_0)^\alpha (1+\delta_1)^{-(1-\alpha)}<e^{-2c}<1.
\end{equation}
We take $\alpha>\gamma_0$.
By hypothesis (H3), $m_2 < M_1$ and $m_2-m_1<\beta$. So,
taking $\beta$ and $\delta_0$ small, and $\alpha$ close
enough to $1$, we ensure that
\begin{equation}
\label{eq.beta}
\alpha m_2 + (1-\alpha)M_2 <  \gamma_0 m_1 + (1-\gamma_0)M_1 -
l \log (1+\delta_0).
\end{equation}

\smallskip

Preparing the definition of $\SK$, we introduce the compact convex set
$K_\alpha\subset\SI$ given by
$$
K_{\alpha}=\{\mu \in \SI; \mu(V)\leq \alpha \}.
$$
Since Lebesgue almost every point spends at most a fraction
$\gamma_0<\alpha$ of time inside $V\subset B_{p+1}\cup \dots \cup B_{p+q}$,
all the ergodic absolutely continuous invariant measures constructed in
\cite{ABV00} belong to $K_\alpha$. In particular, $K_\alpha$ is
non-empty. We will see that $K_\alpha$ contains the equilibrium
states of potentials with low variation.

Let us recall the ergodic decomposition theorem, as it is proven
in \cite{Man87}:

\begin{theorem}\label{EDT}
Given any invariant measure, there are ergodic invariant measures
$\{\mu_x: x\in M\}$ depending measurably on the point $x$ and such
that
$$
\int f d\mu = \int (\int fd\mu_x) d\mu \quad\text{for every $f\in
L^1(d\mu)$.}
$$
Moreover, this decomposition $\{\mu_x: x\in M\}$ is essentially
unique and, in fact, $\mu_x=\lim n^{-1}\sum_{j=0}^{n-1}\delta_{f^j(x)}$
at $\mu$-almost every point.
\end{theorem}

Our distinguished set of invariant measures is the non-empty
 set $\SK\subset K_\alpha$ defined by
\begin{equation*}
\SK = \{ \mu \in \SI; \mu_x \in K_\alpha \text{ for $\mu$-a.e. } x \}
\end{equation*}

\begin{definition}
We say that a measure $\nu$ is \emph{$f$-expanding with exponent
$c$} if for $\nu$-almost every $x\in M$ we have:
$$
\limsup\limits_{n\rightarrow +\infty}
\frac{1}{n}\sum\limits_{j=0}^{n-1} \log \|Df(f^j(x))^{-1}\| \leq
-2c <0.
$$
\end{definition}

\begin{remark} \label{rem1}
If $\nu$ is an ergodic invariant measure, by the ergodic theorem
the above definition is equivalent to $\int \log\| Df^{-1}\|
d\nu \leq -2c.$
\end{remark}

\begin{remark}
Given any ergodic invariant measure $\nu$, the Lyapunov exponents
of $(f,\nu)$ are all positive if and only if $\nu$ is
$f^N$-expanding for some iterate $N$. See \cite{ABV00}. Also, if
every invariant measure $\nu$ is $f$-expanding, then $f$ is a
uniformly expanding map. See \cite{AAS03}. In fact, the same
conclusion holds, more generally, if all invariant measures have
only positive Lyapunov exponents. See \cite{Cao}.
\end{remark}

The next statement proves that all measures in $\SK$ are
$f$-expanding, with uniform exponent. Let $c>0$ be as in
\eqref{eq.alpha}.

\begin{lemma}\label{S1}
Every measure in $\mu \in \SK$ is $f$-expanding with exponent $c$:
$$
\limsup\limits_{n\rightarrow +\infty}
\frac{1}{n}\sum\limits_{j=0}^{n-1} \log \|Df(f^j(x))^{-1}\| \leq
-2c <0 \quad\text{for $\mu$-almost every $x\in M$.}
$$
\end{lemma}

\begin{proof}
Suppose that $\mu \in \SK$ is ergodic. Since $\mu \in K_\alpha$ and
$\mu(V)\leq \alpha$, by the ergodic theorem almost every point
spends a fraction less than $\alpha$ inside $V$: there exist an invariant
set $A$ with $\mu(A)=1$ such that for every $x\in A$,
$$
\lim\limits_{n\rightarrow \infty} \frac{1}{n} \sum\limits_{i=0}^{n-1}
\chi_V(f^i(x)) \leq \alpha.
$$
By hypothesis (H1), we have that  $\|Df(y)^{-1}\| \leq (1+\delta_0)$ for
every $y\in V$.
Moreover, $\|Df(y)^{-1}\|\leq (1+\delta_1)^{-1}$ for $y\in V^c$.
This implies that

\begin{equation*}\begin{aligned}
\frac{1}{n} \sum\limits_{i=0}^{n-1} \log \|Df(f^i(x))^{-1}\|
& \leq  \frac{1}{n}\log (1+\delta_0)^{\alpha n}(1+\delta_1)^{-(1-\alpha) n}
\\
& \le \log(1+\delta_0)^{\alpha}(1+\delta_1)^{-(1-\alpha)} < - 2c < 0
\end{aligned}\end{equation*}
for all $x\in A$.

 To finish , let $H$ be the set of $x$ that satisfy the condition in the
conclusion of the lemma.  Since every $\mu\in \SK$ is $\mu=\{\mu_x\}$ a convex
combination of ergodic measures
$\mu_x$ in $K_\alpha$, by the previous case, $\mu_a(H)=1$ for $\mu$ almost
every
$a$, and this implies
that $\mu(H)=\int \mu_x(H)d\mu=1$.
\end{proof}

Now we need the notion of hyperbolic time, first introduced by
Alves~\cite{Al00}.

\begin{definition}
We say that $n$ is a \textit{hyperbolic time} for $x$ with exponent $c$, if
for
every  $1\leq j\leq n$:
$$
\prod _{k=0}^{j-1}  \|Df(f^{n-k}(x))^{-1}\| \leq e^{-cj}.
$$
\end{definition}

To prove that for an $f$-expanding measure almost every point
admits infinitely many hyperbolic times, we need the following
lemma due to Pliss. See for instance \cite{ABV00} for a proof.

\begin{lemma}
Given $A\geq c_2>c_1>0$, let $d_0=\frac{c_2-c_1}{A-c_1}$. If $a_1,
\dots, a_n$ are real numbers such that $a_i\leq A$ and
$$
\sum\limits_{i=1}^{n} a_i \geq c_2 n
$$
then there are integer numbers $l>d_0 n$ and $1<n_1<\dots<n_l\leq
n$ so that, for every $0\leq k \leq n_i$ and $i=1,\dots,l$ :
$$
\sum\limits_{j=k+1}^{n_i} a_j \geq c_1(n_i-n)
$$
\end{lemma}

Using this lemma, we get

\begin{lemma}
\label{lematemposhiperbolicos} For every invariant measure $\nu$
which is $f$-expanding with exponent $c$, there exists a full
$\nu$-measure set $H\subset M$ such that every $x\in H$ has
infinitely many hyperbolic times $n_i=n_i(x)$ with exponent $c$ and, in
fact, the
density of hyperbolic times at infinity is larger than some
$d_0=d_0(c)>0$:

\begin{enumerate}
\item $\displaystyle{\prod _{k=0}^{j-1}  \|Df^{-1}(f^{n_{i}-k}(x))\| \leq
e^{-cj}}$ for every $0\leq j\leq n_i$
\item $\displaystyle{\liminf\limits_{n\rightarrow \infty}
\frac{\sharp\{0\leq n_i\leq n\}}{n}\geq d_0>0}$.
\end{enumerate}
\end{lemma}

\begin{proof}
By the definition of $f$-expanding measures, there exists a set
$H$ with $\nu(H)=1$ such that given any $x\in H$ we have

$$
\sum\limits_{i=0}^{n-1} \log \|Df(f^i(x))^{-1}\| \leq -\frac{3c}{2} n
$$
for every $n$ large enough. Then, it suffices to take
$A=\sup\limits_{x\in M} -\log \|Df(x)^{-1}\|$, $c_1 = c$, $c_2=\frac{3c}{2}$
and $a_i=-\log \|Df(f^{i-1}(x))^{-1}\|$ in the previous lemma.
\end{proof}

The next lemma asserts that points at hyperbolic times have
unstable manifolds with size uniformly bounded from below.

\begin{lemma} \label{epsilon}
There exists $\epsilon_0>0$ such that for every $x$ and $n_i$ a
hyperbolic time of $x$, if  $z\in M$ satisfies $f^{n_i}(z)\in
B_{\epsilon_0} (f^{n_i}(x))$ then
$$
d(f^{n_i-j}(x),f^{n_i-j}(z))\leq e^{\frac{-c}{2}j}
d(f^{n_i}(z),f^{n_i}(x))
\quad\text{for every $0\le j \le n_i(x)$.}
$$
\end{lemma}

\begin{proof}
Since $Df$ is uniformly continuous and a local diffeomorphism, there exists
$\epsilon_0$ such that for every $\xi,\eta \in M$ with $\xi\in
B_{\epsilon_0}(\eta)$ then
$$
\frac{\|Df(\xi)^{-1}\|}{\|Df(\eta)^{-1}\|}
\leq e^{\frac{c}{2}}.
$$
By definition, if $n_i$ is hyperbolic time for $x$ then $\prod _{k=0}^{j-1}
\|Df(f^{n_{i}-k}(x))^{-1}\| \leq e^{-cj}$ for every $0\leq j\leq n_i.$
Observe that $d(f^{n_i -1}(z),f^{n_i -1}(x))\leq \epsilon_0$.
This is because $f^{n_i}(z) \in B_{\epsilon_0} (f^{n_i}(x))$ and, by
the previous observation, the norm of the derivarive of the inverse
branch of $f$ that sends $f^{n_i}(x)$ to $f^{n_i-1}(x)$ is less
than $e^{-\frac{c}{2}}$ restricted to
$B_{\epsilon_0}(f^{n_i}(x))$. Arguing by
induction,
$$
\prod _{k=0}^{j-1}  \|Df(f^{n_{i}-k}(z))^{-1}\| \leq
e^{\frac{-c}{2}j} \quad\text{for all}\quad 0\leq j\leq n_i.
$$
We conclude that $d(f^{n_i-j}(x),f^{n_i-j}(z))\leq
e^{\frac{-c}{2}j}d(f^{n_i}(x),f^{n_i}(z)),$ proving the lemma.
\end{proof}

Since $c$ is fixed, we will write simply $\nu$-expanding to mean
$\nu$-expanding with exponent $c$.

\section{Existence of equilibrium states for continuous low variation
potentials}

Our main aim in this section is to establish the existence of
equilibrium states for continuous low variation potentials, in order to prove
theorem \ref{theoremA}. In particular, this applies to $\phi=0$,
which always has low variation. Thus, our construction also yields
maximal entropy measures for these transformations.

Beforehand, we use some results of the previous section to establish
expansiveness for measures in $\SK$.

\begin{definition}
Given $\epsilon>0$ we define the set $A_\epsilon(x)$ by:

$$
A_\epsilon(x)=\{y\in M;d(f^n(x),f^n(y))\leq \epsilon \text{ for every }
n\geq 0\}.
$$
\end{definition}

By definition, $f$ is an \textbf{expansive map} with expansiveness
constant $\tilde{\epsilon}$ if and only if
$A_\epsilon(x)=\{x\}$ for every $x\in M$ and
$\epsilon<\tilde{\epsilon}$.

\begin{lemma}\label{l2}
Suppose that $\mu \in \SK$  and let $\epsilon_0$ be as
constructed in lemma \ref{epsilon}. Then for $\mu$-almost every
$x\in M$ and any $\epsilon<\epsilon_0$,
$$
A_\epsilon(x)=\{x\}.
$$
\end{lemma}

\begin{proof}
By lemma \ref{lematemposhiperbolicos} we have that almost
every $x\in M$ has infinitely many hyperbolic times $n_i(x)$.
By lemma \ref{epsilon}, if $z\in A_\epsilon(x)$ with
$\epsilon<\epsilon_0$ then for any $n_i$ we have
$$
d(x,z)\leq e^{-\frac{c}{2} n_i}d(f^{n_i}(x),f^{n_i}(z))
\leq e^{-\frac{c}{2} n_i}\epsilon.
$$
Making $n_i \rightarrow \infty$ we deduce that $x=z$.
\end{proof}

Let $\SP=\{P_1,\dots,P_l\}$ be any partition of $M$ in measurable
sets with diameter less than $\epsilon_0$. From the above lemma we
get

\begin{corollary}\label{c1}
$\SP$ is a generating partition for every $\mu \in \SK$.
\end{corollary}

\begin{proof} Define
$$
\SP^{(n)}=\{C^{(n)}=P_{i_0}\cap \dots \cap f^{-n+1}(P_{i_{n-1}})\},
\quad\text{for each $n\ge 1$.}
$$
We need to prove that given any measurable set $A$ and given
$\delta>0$ there exist elements $C_1^{(n)}, \dots,
C_m^{(n)}$ of $\SP^{(n)}$ such that
$$
\mu(\bigcup C_i^{(n)}\Delta A)\leq \delta.
$$
Consider  $K_1\subset A$ and $K_2\subset A^c$  compact sets  such that
$\mu(K_1 \Delta A) \leq \delta$ and $\mu(K_2 \Delta A^c) \leq \delta$. Let
$r=d(K_1,K_2)>0$. Lemma \ref{l2} gives that if $n$ is big enough then
$\text{diam}\SP^{(n)}(x)\leq \frac{r}{2}$, for $x$ in a set with
$\mu$-measure bigger than $1-\delta$.
Consider the sets $C_1^{(n)}, \dots, C_m^{(n)} \in \SP^{(n)}$ that intersect
$K_1$.
Then
\begin{equation*}\begin{aligned}
\mu( \bigcup C_i^{(n)}\Delta A)
& = \mu(\bigcup C_i^{(n)}-A) + \mu(A - \bigcup C_i^{(n)})
\\ & \leq \mu(A-K_1)+\mu(A^c-K_2)
  + \delta \leq 3\delta.
\end{aligned}\end{equation*}
This proves the claim.
\end{proof}

\begin{remark}
\label{rem.us}
Recalling the definition of entropy with respect to a partition
$\SQ$,
$$
H_\mu(\SQ)=\sum\limits_{Q\in \SQ} -\mu(Q)\log \mu(Q),
$$
we have that for any partition $\SQ$ such $\mu_0(\partial Q)=0$, for each $Q\in
\SQ$, then the function $\mu \rightarrow
H_\mu(\SQ)$ is continuous in $\mu_0$. This imply that

$$
\mu \mapsto h_\mu(f,\SP)=\inf\limits_{n\rightarrow \infty}
    \frac{1}{n}H_\mu(\SP^{(n)}).
$$
is upper-semicontinuous in $\mu_0$.
\end{remark}

Observe that, by corollary \ref{c1}, $\SP$ is a
generating partition for every $\mu\in\SK$. So, as a consequence of
 Kolmogorov-Sinai's theorem(see e.g. \cite{Man87}),

\begin{corollary}\label{c2}
For every measure $\mu \in \SK$ we have $h_\mu(f)=h_\mu(f, \SP)$.
\end{corollary}

The next lemma  will allow us prove that if $\phi$ has $\rho$-low
variation, the function  $F_\phi$ restrict to $\SK$ has a maximum $\mu_\phi$
and
this maximum
is in fact an equilibrium state for $\phi$.

\begin{lemma}\label{S2}
All ergodic measures $\eta$ outside $\SK$ have small entropy:
there exists $\rho<1$ such that if $\eta \in \SK^c$ is ergodic
then
$$
h_\eta (f) \leq \rho h_{top}(f).
$$
\end{lemma}

\begin{proof}

As we are supposing that $\eta$ is ergodic, we have that
$\eta(V)>\alpha$, since every ergodic measure $\mu$ such
$\mu(V)\leq \alpha$ lies in $\SK$. Denoting $\lambda_1(x)\geq
\lambda_2(x) \geq \dots\lambda_s \geq 0>\lambda_{s+1}\dots \geq
\lambda_l(x)$ the Lyapunov exponents in $x$, we know that
$\lambda_i=\lambda_i(x)$ is constant $\eta$-almost everywhere. By
the theorem of Oseledets~\cite{Os68},
\begin{equation}
\label{eq.star1}
\int \log\|\det Df\|d\eta = \sum\limits_{i=1}^l \lambda_i.
\end{equation}
On the other hand, we have that $\lambda_l>-\log(1+\delta_0)$,
since by hyphotesis, $\|Df(x)^{-1}\|\leq 1+\delta_0$. By Ruelle's
inequality (see \cite{Ru78}) we have that
\begin{equation}
\label{eq.star2}
h_\eta(f) \leq \sum\limits_{i=1}^s \lambda_i= \int \log \|\det Df(x)\|d\eta -
\sum\limits_{i=s+1}^l
\lambda_i.
\end{equation}
Since $m_2=\sup\limits_{x\in V}  \log \|\det Df(x)\| <
M_2=\sup\limits_{x\in V^c} \log \|\det Df(x)\|$  and
$\eta(V)>\alpha$  we have:
\begin{equation*}\begin{aligned}
h_\eta(f) \leq \int \log|\det Df|d\eta
 & \leq \eta(V)m_2 + (1-\eta(V))M_2 + (l-s)\log(1+\delta_0) \\
 & \leq \alpha m_2 + (1-\alpha)M_2 +l\log(1+\delta_0)
\end{aligned}\end{equation*}
Let $\mu_0$ be any ergodic absolutely continuous invariant measure
as constructed in \cite{ABV00}. Since Lebesgue almost every point
spends at most a fraction $\gamma_0$ of time inside $W\subset
B_{p+1}\cup \dots \cup B_{p+q}$, we have that $\mu_0(W)<\gamma_0$.
As $f$ is $C^{1+\alpha}$ and $\mu_0$ is absolutely continuous, we
may use Pesin's entropy formula (see \cite{Pe77}):
$$
h_{\mu_0}(f)= \int \log \|\det Df\|d\mu_0 \geq \mu_0(W)m_1 + (1-\mu_0(W))M_1
$$
As $\mu_0(W)\leq \gamma_0 $ and $m_1<M_1$, we conclude that
$$
\gamma_0 m_1 + (1-\gamma_0)M_1 \leq  h_{\mu_0}(f).
$$
By  \eqref{eq.beta},
$$
\alpha m_2 + (1-\alpha)M_2 < \gamma_0 m_1 + (1-\gamma_0)M_1 -
l\log(1+\delta_0).
$$
So, we can choose $\rho<1$ such that
\begin{equation*}\begin{aligned}
h_\eta(f)
 & \leq \alpha m_2 + (1-\alpha)M_2 + l\log(1+\delta_0)
 < \rho(\gamma_0 m_1 + (1-\gamma_0)M_1 ) \\
 & < \rho h_{\mu_0}(f) \leq \rho h_{top}(f).
\end{aligned}\end{equation*}
This proves lemma \ref{S2}.
\end{proof}

\begin{remark} Observe that if $\mu_0$ is some SRB measure for $f$, from the
proof of previous lemma, we may choose  $\rho <
\frac{h_{\mu_0}}{h_{top}(f)}$.
\end{remark} 

Observe that it follows  from the lemma \ref{S2} and the variational principle:

\begin{corollary}[Variational Principle for expanding measures]
 \label{VP}
If $\phi$ is a $\rho$-low variation potential, then:

$$
\sup\limits_{\nu\in \SK} h_\nu(f) + \int \phi d\nu =P(\phi)
$$

In particular,

$$\sup\limits_{\nu\in \SK} h_\nu(f) = h_{top}(f)$$

\end{corollary}

\begin{proof}

  Denote by $E$ the set of all ergodic invariant probabilities, to prove the
lemma, we just need to prove  that:

  \begin{equation}\label{eq.er}
  \sup\limits_{\nu\in \SK} F_\phi(\nu) =
  \sup\limits_{\nu\in  E} F_\phi(\nu)
  \end{equation}

 Since

  $$P(\phi) =
  \sup\limits_{\nu\in  E} F_\phi(\nu) $$

To prove \ref{eq.er}, note that  by lemma \ref{S2} we have that for $\nu\in
\SK^c$  ergodic, then $h_{\nu}(f)\leq \rho h_{top}(f)$. This imply that:

$$
F_{\phi}(\nu)=h_{\nu}(f) + \int \phi d\nu\leq \rho h_{top}(f) + \max_{x\in M}
\phi(x)  < P(\phi)
$$

 and this prove the corollary.

\end{proof}

\section{Proof of theorem \ref{theoremA}}

First of all, observe that the previous corollary asserts that
$\sup\limits_{\nu\in \SK} F_\phi(\nu)= P(\phi)$. To prove that there exists
some
equilibrium states,  consider a sequence of measures $\mu_k \in \SK$ such
$F_\phi(\mu_k)$ converge to $P(\phi)$.

Without loss of generality, we suppose that $\mu_k \rightarrow \mu$ weakly.  We
prove that $\mu$ is an equilibrium state for $\phi$ and belongs to $\SK$.

First, we claim that $F_\phi(\mu)= P(\phi)$. In fact, fix $\SP$ a
 partition with diameter less than $\epsilon_0$,  such $\mu(\partial P)=0$ for
any $P\in \SP$. Observe that, since $\mu_k\in \SK$, we have that
$h_{\mu_k}(f)=h_{\mu_k}(f,\SP)$, by corollary \ref{c2}. Then:

$$
P(\phi)=\sup\limits_{\nu\in \SK} F_\phi(\nu) = \limsup F_\phi(\mu_k) =
\limsup h_{\mu_k}(f) +\int \phi d\mu_k $$

As, by remark \ref{rem.us}, the function $\nu \rightarrow h_\nu(f,\SP)$
is upper-semicontinuous in $\mu$, thus:

$$
 \limsup h_{\mu_k}(f,\SP) + \int\phi d\mu_k
\leq h_\mu(f,\SP) + \int \phi d\mu \leq h_\mu(f) + \int \phi d\mu = F_\phi(\mu)
$$

Then, we have

$$P(\phi)=\sup\limits_{\nu\in \SK} F_\phi(\nu)\leq F_\phi(\mu)\leq P(\phi),$$

which imply that $\mu$ is equilibrium state for $\phi$.

Now, we prove that any  measure $\eta$  such $F_\phi(\eta)= \sup\limits_{\nu\in
\SK}
F_\phi(\nu)$ belongs to $\SK$, proving that all equilibrium states belong to
$\SK$.

  In fact, if $\eta=\{\eta_x\}$ is the ergodic decomposition of $\eta$, we
should prove that the set  $S=\{x\in M; \eta_x \in K_\alpha\}$ is a $\eta$-full
measure set.

  Using the fact that

$$
h_\eta(f)=\int h_{\eta_x}(f)\,d\eta(x),
$$
(see \cite{Ro67}, for instance), we have $F_\phi(\eta) = \int F_\phi(\eta_x)\,
d\eta(x)$

  Suppose by contradiction
  that $\eta(S^c)>0$. Observe that if $y \in S^c$, then $\eta_y$ is in the
  hyphoteses of lemma \ref{S2} and thus:

 \begin{equation}\label{eq.max}
F_\phi(\eta_y) = h_{\eta_y}+\int\phi\,d\eta_y
< \rho h_{top}(f) + \max_{x\in M} \phi(x)  <P(\phi),
  \end{equation}

   since $\phi$ is a $\rho$-low variation  potential.

  As for every $x\in S$ we have that $F_\phi(\mu_x)\leq P(\phi)$, the
inequality
\ref{eq.max} implies that $F_\phi(\eta) = \int F_\phi(\eta_x)\,
d\eta(x)<P(\phi)$, which is a contradiction.
   Then, $\eta(S)=1$ which imply, by the definition of $\SK$ that $\eta\in \SK$
and that all equilibrium states are in $\SK$, completing the prove of theorem
\ref{theoremA}.

\section{Proof of Theorem B  }

We give a proof of the existence
and uniqueness of a measure with maximum entropy. Throughout, we assume
the additional hypothesis (H4): existence of a transitive Markov partition.
Observe that we do not require this partition to be generating.

Firstly, by Theorem~\ref{theoremA}, there exists some measure $\mu_{max}$ with
maximal entropy. Our strategy to prove its unicity is  transfer our problem to a
subshift of finite type $\sigma^+:\Sigma^+\to\Sigma^+$ via ergodic conjugacy. If
we denote $\partial R = \bigcup\limits_{n=1}^d \partial R_i$ and
$\tilde{M}=M-\bigcup\limits_{n\geq 0}f^{-n}(\partial R)$, we have that if $\mu$
is an ergodic measure such $\mu(\partial \SR)=1$ then the entropy of $\mu$ is
less than topological entropy, since $h_{top}(f|\partial \SR)<h_{top}(f)$.
Thus, we just consider invariant probabilities such $\mu(\tilde{M})=1$. We may
define a map $\Pi: \tilde{M}\rightarrow \Sigma^+$
over a subshift of finite type $\Sigma^+$ associate to some transition matrix
$A$, by
$$
\Pi(x)=(i_0,\dots,i_n,\dots) \text{ such that } f^n(x)\in P_{i_n}.
$$
Observe that this map is  a semiconjugacy between $f$
and $\sigma^+$.    Define the cylinders
$$[i]=[i_0, \dots, i_n]=\{x\in M; 0\leq j\leq n,
f^j(x)\in R_{i_j}\}.
$$

\begin{definition}
Let $(i_n)$ be the itinerary of $x$, defined by $f^n(x)\in
R_{i_n}$, for each $n\geq 0$. We define $[x]$ to be the set
$$
[x]=[i_0,...,i_n,...]= \{y\in M; f^n(y)\in R_{i_n}\}.
$$
\end{definition}

Given an invariant measure $\eta$ satisfying $[x]=\{x\}$ for $\eta$-a.e., then
$\Pi$ is an ergodic conjugacy between $(f,\eta)$ and
$(\sigma^+,\Pi^{\star}\eta)$, where $\Pi^{\star}\eta)$ is defined by
$\Pi^{\star}\eta)(A) = \eta(\Pi^{-1}(A))$.
Observe that  some measures can not be transported to the shift but, by lemma
\ref{S2}  any $f$-invariant measure with big entropy has hyperbolic times for
almost everywhere. It allow us to prove $[x]=\{x\}$ for $\eta$-a.e., which imply
that $\eta$ is ergodically equivalent to some measure in the shift.

Using the classical fact that transitive subshifts of finite
type have exactly one measure of maximal entropy, we prove that $f$ admits only
one measure $\mu_{max}$ with maximal entropy. If $\Sigma^+$ is topologically
mixing, then its maximal measure is mixing (\cite{Bo75}). Since, by Ornstein's
Theorem(\cite{Man87}), every mixing Markov measure is Bernoulli, we have that
$\mu_{max}$ is Bernoulli.

\vspace{1cm}

\noindent Krerley Oliveira ( krerley{\@@}mat.ufal.br )\\
Departamento de Matematica - UFAL, Campus A.C. Simoes, s/n 57072-090 Maceio,
Alagoas
- Brazil

\end{document}